\documentclass{amsart}

\usepackage{amsfonts}
\usepackage{amsthm,amsmath}
\usepackage{amsbsy}
\usepackage{bm}
\usepackage{amssymb} 
\usepackage{verbatim}
\RequirePackage[dvips]{hyperref}
 \usepackage{lineno}
 \usepackage{color}

\usepackage[mathscr]{eucal}
\usepackage{times}

\newcommand{\E}{\mbox{\bf E}}
\def\P{{\bf P}}
\def\Var{\mbox{Var}}
\newcommand{\Cov}{\mbox{\rm Cov}}

\def\half{\frac{1}{2}}

\newcommand{\one}{{\mathbf 1}}

\def\intt{\int\limits}

\def\prodd{\prod\limits}

\def\to{\rightarrow}

\def\mb{\mbox}

\newcommand{\para}[1]{\vspace{4mm}\noindent{\bf #1:}}

\def\l{\left}
\def\r{\right}
\def\<{\langle}
\def\>{\rangle}

\newcommand{\ba}{\[\begin{aligned}}
\newcommand{\ea}{\end{aligned}\]}
\newcommand\mnote[1]{} 
\newcommand{\beq}[1]{\begin{equation}\label{#1}}
\newcommand\eeq{\end{equation}}
\newcommand\ben{\begin{equation}}
\newcommand\een{\end{equation}}
\newcommand\besn{\begin{eqnarray}}
\newcommand\eesn{\end{eqnarray}}

\def\bthm{\begin{theorem}}
\def\ethm{\end{theorem}}
\def\bdefn{\begin{definition}}
\def\edefn{\end{definition}}
\def\benu{\begin{enumerate}}
\def\eenu{\end{enumerate}}
\def\beit{\begin{itemize}}
\def\eeit{\end{itemize}}
\def\beds{\begin{description}}
\def\eeds{\end{description}}
\def\bepr{\begin{problem}}
\def\eepr{\end{problem}}

\def\bprf{\begin{proof}}
\def\eprf{\end{proof}}
\def\berk{\begin{remark}}
\def\eerk{\end{remark}}
\def\bex{\begin{exercise}}
\def\eex{\end{exercise}}
\def\beg{\begin{example}}
\def\eeg{\end{example}}

\def\suchthat{{\; : \;}}

\def\R{\mathbb{R}}

\def\Z{\mathbb{Z}}

\renewcommand{\phi}{\varphi}

\def\bet{\beta}
\def\gam{\gamma}
\def\del{\delta}
\def\eps{\epsilon}

\def\lam{\lambda}

\def\sig{\sigma}

\theoremstyle{plain} 
    \newtheorem{theorem}{Theorem}
    \newtheorem{lemma}[theorem]{Lemma}

\theoremstyle{definition} 
    \newtheorem{definition}[theorem]{Definition}

    \newtheorem{remark}[theorem]{Remark}
    \newtheorem{example}[theorem]{Example}
    
    \newtheorem{exercise}[theorem]{Exercise}
   \newtheorem{problem}[theorem]{\bf }

\def\half{\frac{1}{2}}

\begin{document}

\title[Persistence probabilities in Gaussian processes]{Persistence probabilities in centered, stationary, Gaussian processes in discrete time}

\author{Krishna M. and  Manjunath Krishnapur}
\address{The Institute of Mathematical Sciences \\
Taramani, Chennai 600113, India}

\email{krishna@imsc.res.in}

\address{Department of Mathematics\\
        Indian Institute of Science\\
        Bangalore 560012, India}

\email{manju@math.iisc.ernet.in}

\thanks{Partially supported by IMSc Project 12-R\&D-IMS-5.01-0106 and UGC center for advanced studies.
  }

\def\T{{\mathbb T}}

\date{\today}


\maketitle

\section{The problem and  our results}
Let $X=(X_{m})_{m\in \Z^{d}}$ be a centered, stationary Gaussian process on $\Z^{d}$. This means that for any $k\ge 1$ and any $m_{0},\ldots ,m_{k}\in \Z^{d}$, the vector $(X_{m_{j}+m_{0}})_{1\le j\le k}$ has a multivariate Gaussian distribution with zero mean and a covariance matrix that does not depend on $m_{0}$. For basics on Gaussian processes, consult for example, the book by Adler~\cite{adler}.

For a subset $A\subseteq\Z^{d}$, we define the {\em persistence probability} (also called gap probability or hole probability) of $X$ in $A$ as
\ba
H_{X}(A):=\P\{X_{m}>0\mb{ for all }m\in A\}.
\ea
In particular, one may be interested in $H_{X}(Q_{N})$, where for $N=(N_{1},\ldots ,N_{d})\in \Z^{d}$, the cube $Q_{N}:=\{m\in \Z^{d}\suchthat 1\le m_{k}\le N_{k}\mb{ for each }1\le k\le d\}$.   This paper is exclusively about getting bounds on the persistence probability under some additional conditions on the Gaussian process.

\para{Notations}  Let $\T=[-\pi,\pi]$. For $d\ge 1$, we use $\lam$  to denote the Lebesgue measure on $\T^{d}$ normalized so that $\lam(\T^{d})=1$. A stationary Gaussian process on $\Z^{d}$ is uniquely described by its covariance kernel $\Cov(X_{m},X_{n})=\Cov(X_{0},X_{m-n})$.  Further, there exists a unique finite Borel measure $\mu$ on $\T^{d}$ that is symmetric about the origin (i.e., $\mu(I)=\mu(-I)$ for any Borel set $I\subseteq \T^{d}$) such that  $\Cov(X_{0},X_{m})=\hat{\mu}(m)$ where $\hat{\mu}(m)=\int_{\T^{d}}e^{i\<m,t\>}d\mu(t)$ with the usual notation for the inner product $\<m,t\>=m_{1}t_{1}+\ldots +m_{d}t_{d}$. The measure $\mu$ is called the spectral measure of the process $X$. Write $d\mu(t)=b(t)d\lam(t)+d\mu_{s}(t)$ where $\mu_{s}$ is singular to Lebesgue measure and $b\in L^{1}(\T^{d},\lam)$ is non-negative. In all the results of this paper, it will be assumed that $b$ is not identically zero. In other words, the spectral measure is not singular. Lastly, for a subset $A\subseteq \Z^{d}$, we denote the covariance matrix of $(X_{m})_{m\in A}$ by $\Sigma_{A}:=\l(\hat{\mu}(j-k)\r)_{j,k\in A}$ and the cardinality of $A$ by $|A|$.  

These notations will be maintained throughout the paper without further mention. In addition, there will appear many constants denoted by $C,c,\gam$ etc. Unless otherwise mentioned, the constants  depend on the given process $X$ (or equivalently, on the spectral measure $\mu$).

We now state our results and then give an overview of past results in the literature in Section~\ref{sec:reviewofresults}. Our first theorem has already been proved by N. Feldheim and O. Feldheim~\cite{feldheims} and but we explain in Section~\ref{sec:reviewofresults} why we include it here nevertheless. 

\begin{theorem}\label{thm:lbdexponential} Assume that $b(t)\ge \del$ for $a.e.$ $t\in [-\eps,\eps]^{d}$ for some positive numbers $\del,\eps$. Then, for any finite $A\subseteq \Z^{d}$, we have $H_{X}(A)\ge e^{-\gam |A|}$ for some finite constant $\gam$ that depends only on $\del$ and $\eps$.
\end{theorem}
Theorem~\ref{thm:lbdexponential} is proved in Section~\ref{sec:thmlbdexp}. Then, in Section~\ref{sec:counterexample1}, we exhibit a Gaussian process on $\Z$ that does not satisfy the conditions of this theorem  and for which $H_{X}(\{1,2,\ldots ,N\})$ decays faster than exponentially in $N$. To deal with such cases, in one dimension, we prove different lower bounds under weaker assumptions on the spectral measure. With a slight abuse of notation, we write  $H_{X}(N)$ for $H_{X}(\{1,2,\ldots ,N\})$.
\begin{theorem}\label{thm:lbdsubexponential} Let $d=1$. 
\benu
\item If $b(\cdot)$ is not identically zero in $L^{1}(\T)$, then $H_{X}(N)\ge e^{-\gam N^{2}}$ for all $N$, for some finite constant $\gam$ that may depend on $b(\cdot)$. 
\item Assume that there exist some $p>0$ and $C<\infty$ such that $\lam\{t\in \T^{d}\suchthat b(t)\le \del\}\le C\del^{p}$ for any $\del\in [0,\pi]$. Then, $H_{X}(N)\ge e^{-\gam N\log N}$ for all $N$, for some finite constant $\gam$ that depends only on $C$ and $p$.
\eenu
\end{theorem}
Theorem~\ref{thm:lbdsubexponential} is proved in Section~\ref{sec:lbdsubexponential}.
Lastly, in Section~\ref{sec:counterexample2}, we give an example of a Gaussian process on $\Z$ for which the gap probability appears to achieve the lower bound of $e^{-\gam N^{2}}$. We do not have a full proof that it works, but we present convincing numerical evidence.

\section{Brief review of past results}\label{sec:reviewofresults}
 As may be expected, the question of gap probability has been studied quite extensively. We give a brief overview of some of the relevant results and then explain where our results fit in. 
 
 Early papers on gap probabilities are by  Longuet-Higgins \cite{LoHi} and Newell and Rosenblatt ~ \cite{newellrosenblatt}.  Newell and Rosenblatt~\cite{newellrosenblatt}  obtain a number of bounds for the gap probability $\P\{X_{t}>0\mb{ for }t\in [0,T]\}$ for a stationary Gaussian process on $\R$. They showed that if the covariance of $X_{0}$ and $X_{m}$ goes to zero as $m\to \infty$, then  the persistence probability $H_{X}(Q_{N})$ decays faster than any polynomial in $N$ and if $\Cov(X_{0},X_{m})$ is also summable, then they showed that $H_{X}(N)\le e^{-cN}$. They also obtained lower bounds, but under the assumption that the covariance is positive.  Some of these results generalize to higher dimensions, see for example the paper of Malevich~\cite{malevich}. A more recent paper of Dembo and Mukherjee~\cite{dembomukherjee} is also concerned with the question of gap probability for one-dimensional Gaussian processes,  but again they assume positivity of covariance. In all these papers, the assumption of positive covariance is crucial in that it allows one to compare with other processes (for example, the i.i.d. process)  using Slepian's inequality (this inequality is recalled in Section~\ref{sec:thmlbdexp}).
 
Our interest is in getting lower bounds for the gap probability even when the covariance is not positive. The first result we know of this kind is due to Antezana, Buckley, Marzo and Olsen~\cite{abmo} who showed that for the Gaussian process $(X_{t})_{t\in \R}$ with $\Cov(X_{t},X_{s})=\frac{\sin(t-s)}{t-s}$ (known as the Paley-Wiener process), the gap probability has the bounds $e^{-c_{1}T}\le \P\{X_{t}>0\mb{ for }0\le t\le T\}\le e^{-c_{2}T}$. 

Generalizing their result, N. Feldheim and O. Feldheim~\cite{feldheims} showed similar exponential upper and lower bounds for a large class of Gaussian processes in $\R$ or $\Z$. Their conditions are similar to ours and their result is stronger than Theorem~\ref{thm:lbdexponential}. The strength is in that lower bounds for gap probability for a process $X=(X_{t})_{t\in \R}$ in continuous time imply also a lower bound for a discrete time process $(X_{n\del})_{n\in \Z}$ obtained by sampling the continuous time process at regular intervals. But it is not possible to go in the reverse direction, since our results will only give $\P\{X_{n\del}>0\mb{ for }n\del\le T\}\ge e^{-c\del^{-1}T}$ and the constant in the exponent blows up when $\del\to 0$. 

Their proof of lower bound uses the result of  \cite{abmo} for the Paley-Wiener process. But we give a full proof of Theorem~\ref{thm:lbdexponential} as it is different and self-contained. We feel that it may have some interesting points to it (in particular, the construction in Lemma~\ref{lem:existenceofh}). 

Hole probabilities have also been studied beyond the setting of stationary Gaussian processes. For example, persistence probability of a random polynomial with i.i.d. coefficients was studied in \cite{dembopoonenshaozeitouni}. Persistence for (the absolute value of) the planar Gaussian analytic function was studied by Sodin and Tsirelson \cite{sotsi} and Nishry\cite{nis}. Also of interest are the results  of Shao and Wang \cite{sha-wan}.

For more on such problems, we refer to the surveys of Li and Shao~\cite{lishao}, the recent review by Frank Aurzada and Thomas Simon \cite{aur-simon} on the persistence question for general L\'evy processes, and the works of Ehrhardt, Majumdar and Bray \cite{bme} and Aurzada and Gullotin-Plantard~\cite{aur-plan} on persistence exponents.

\section{Proof of Theorem~\ref{thm:lbdexponential}}\label{sec:thmlbdexp}
We first present a simple lemma that we shall use many times. The setting and notation are as before. In particular, recall that for $A\subseteq \Z^{d}$, the matrix $\Sigma_{A}:=\l(\hat{\mu}(j-k)\r)_{j,k\in A}$ denotes the covariance matrix of $(X_{m})_{m\in A}$.
\begin{lemma}\label{lem:bboundedbelow} Let $A\subseteq \Z^{d}$ be a finite non-empty subset. Suppose $0\le b_{*},b^{*}\le  +\infty$ are such that $b_{*}\le b(t)\le b^{*}$ for $a.e.$ $t\in \T^{d}$.
\benu
\item Then, $\det(\Sigma_{A})\le \hat{\mu}(0)^{n}$.
\item All the eigenvalues of $\Sigma_{A}$ lie in the interval $[b_{*},\infty)$. If the singular part of the spectral measure $\mu_{s}$ vanishes, then all the eigenvalues of $\Sigma_{A}$ lie in $[b_{*},b^{*}]$.
\item $H_{X}(A)\ge \l(\frac{\sig_{A}}{4\hat{\mu}(0)}\r)^{|A|/2}$ where $\sig_{A}$ is the smallest eigenvalue of $\Sigma_{A}$. In particular, $H_{X}(A)\ge \l(\frac{b_{*}}{4\hat{\mu}(0)}\r)^{|A|/2}$.
\eenu
\end{lemma}
\bprf \benu
\item Since $\Sigma_{A}$ is a positive semidefinite matrix, its determinant is bounded from above by the product of its diagonal entries (we may realise $\Sigma_{A}$ as the Gram matrix of $n$ vectors in $\R^{n}$ and then $\det(\Sigma_{A})$ is the squared volume of the parallelepiped formed by these vectors while the diagonal entries are the squared norms of these vectors). All its diagonal entries are equal to $\hat{\mu}(0)$ and hence the claim follows.
\item Let $u\in \R^{A}$ and set $U(t)=\sum_{k\in A}u_{k}e^{i\<k,t\>}$ for $t\in \T^{d}$. Observe that for $k,\ell\in \Z^{d}$, the inner product of $e^{i\<k,t\>}$ and $e^{i\<\ell,t\>}$ is equal to $\del_{k,\ell}$ in $L^{2}(\T^{d},\lam)$ and equal to $\hat{\mu}(k-\ell)$ in $L^{2}(\T^{d},\mu)$. It  follows that
\ba
\|u\|^{2} &= \int_{\T^{d}}|U(t)|^{2}d\lam(t) \qquad \mb{ and } \qquad
u^{t}\Sigma_{A}u&= \int_{\T^{d}}|U(t)|^{2}d\mu(t).
\ea
Now, $d\mu(t)\ge b(t)d\lam(t)$ with equality if $\mu_{s}=0$. Hence, 
\ba
u^{t}\Sigma_{A}u &\ge \int_{\T^{d}}|U(t)|^{2}b(t)d\lam(t)
\ea
with equality if $\mu_{s}=0$. Thus, using the lower bound for $b$, we see that $u^{t}\Sigma_{A}u\ge b_{*}\|u\|^{2}$. When $\mu_{s}=0$, we may also use the upper bound for $b$ and get $u^{t}\Sigma_{A}u\le b^{*}\|u\|^{2}$. From the variational characterization of eigenvalues of symmetric matrices, the claims follow. 
\item We may assume that $\Sigma_{A}$ is non-singular (otherwise, $\sig_{A}=0$ and by the previous part we must have $b_{*}=0$ and thus the right hand sides of both  inequalities to prove are  zero anyway). Then, the Gaussian vector $(X_{m})_{m\in A}$ has density $(2\pi)^{-|A|/2}\exp\{-\half u^{t}\Sigma_{A}^{-1}u\}$ with respect to Lebesgue measure on $\R^{A}$ and hence
\ba
H_{X}(A) &= \int_{\R_{+}^{d}}\exp\{-\half u^{t}\Sigma_{A}^{-1}u\}\frac{du}{(2\pi)^{|A|/2}\sqrt{\det(\Sigma_{A})}}
\ea
where $\R_{+}=[0,\infty)$. From the first part, we have $\det(\Sigma_{A})\le \hat{\mu}(0)^{|A|}$. Further, $u^{t}\Sigma_{A}^{-1}u\le \frac{1}{\sig_{A}}\|u\|^{2}$ for all $u$ since $\sig_{A}$ is the smallest eigenvalue of $\Sigma_{A}$. Putting these together, we get
\ba
H_{X}(A) &\ge \frac{1}{\hat{\mu}(0)^{|A|/2}}\int_{\R_{+}^{d}}e^{-\frac{1}{2\sig_{A}}\|u\|^{2}}\frac{du}{(2\pi)^{|A|/2}} \\
&= \frac{\sig_{A}^{|A|/2}}{2^{|A|}\hat{\mu}(0)^{|A|/2}}
\ea
by evaluating the integral (which splits into a product of one dimensional Gaussian integrals). This proves the first inequality for $H_{X}(A)$. By the second part, we have the bound $\sig_{A}\ge b_{*}$ from which the second inequality follows. \qedhere
\eenu
\eprf

%


If $b$ is bounded below by a positive constant $b_{*}$ on $\T^{d}$, then the third part of Lemma~\ref{lem:bboundedbelow}  immediately implies the conclusion of Theorem~\ref{thm:lbdexponential} with $\gam=-\half\log(b_{*}/4\hat{\mu}(0))$. But the assumption in that theorem is only that $b$ is bounded below by a positive constant in a neighbourhood of the origin. To get an exponential lower bound under this weaker assumption, we shall use the following direct consequence of Slepian's inequality (see Corollary 2.4 of Adler~\cite{adler} or Slepian's original paper~\cite{slepian}).

\para{Slepian's inequality} Let $X$ and $Y$ be two Gaussian processes on $\Z^{d}$ with $\Cov(X_{m},X_{n})\ge \Cov(Y_{m},Y_{n})$ for all $m,n\in \Z^{d}$ and such that $\Var(X_{m})=\Var(Y_{m})$ for all $m\in \Z^{d}$. Then, $H_{X}(A)\ge H_{Y}(A)$ for any $A\subseteq \Z^{d}$.

\medskip
The idea will be to get a different process $Y$ that is comparable to $X$ as in Slepian's inequality and such that $Y$ has spectral density that is bounded below on all of $\T^{d}$. Then we may combine Slepian's inequality and the lower bound for $H_{Y}(A)$ from Lemma~\ref{lem:bboundedbelow} to prove Theorem~\ref{thm:lbdexponential}. To produce such a $Y$, we shall need the following lemma.

\begin{lemma}\label{lem:existenceofh}
Given $\eps>0$ there exists a function $h\in L^{1}(\T^{d},\lam)$ such that 
\benu
 \item $\hat{h}(n)\ge 0$ for all $n$ and $\hat{h}(0)=0$.
 \item $\sup\limits_{t\in \T^{d}}h(t)=h(0)=1$.
 \item $\sup\limits_{t\not\in (\eps \T)^{d}}h(t)\le -\bet_{\eps}$ where $\beta_{\eps}=\eps^{d}$.
 \eenu
\end{lemma}
Assuming this lemma, we prove Theorem~\ref{thm:lbdexponential}.
\bprf[Proof of Theorem~\ref{thm:lbdexponential}] Let $h$ be the function provided by Lemma~\ref{lem:existenceofh} and set  $\tilde{b}(t)=b(t)-\frac{\del}{2} h(t)$.  Let $d\tilde{\mu}(t)=d\mu_{s}(t)+\tilde{b}(t)dt$  and let $(\tilde{X}_{n})_{n\in \Z}$ be the centered stationary Gaussian process with spectral measure $\tilde{\mu}$.  

From the assumption on $b$, it follows that $\tilde{b}\ge \frac{1}{2}\del$ for  $t\in (\eps \T)^{d}$ and $b(t)\ge \bet_{\eps}$ for  $t\not\in (\eps\T)^{d}$ (for $a.e.$ $t$). Thus, $\tilde{b}$ is bounded below by $\tilde{b}_{*}=\min\{\frac{1}{2}\del, \bet_{\eps}\}$. 

By the properties of $h$, we see that $\hat{\mu}(n)\ge \hat{\tilde{\mu}}(n)$ for all $n$ with equality for $n=0$. In terms of covariances this says that $\E[X_{n}^{2}]=\E[\tilde{X}_{n}^{2}]$ and $\Cov(X_{n},X_{m})\ge \Cov(\tilde{X}_{n},\tilde{X}_{m})$. Therefore, Slepian's inequality applies to $X$ and $\tilde{X}$ and gives $H_{X}(A)\ge H_{\tilde{X}}(A)$.

By Lemma~\ref{lem:bboundedbelow} and the fact that $\hat{\tilde{\mu}}(0)=\hat{\mu}(0)$, we have $H_{\tilde{X}}(A)\ge \l(\frac{\tilde{b}_{*}}{4\hat{\mu}(0)}\r)^{|A|/2}$, completing the proof of the theorem with $\gam=-\half\log(\tilde{b}_{*}/4\hat{\mu}(0))$.  
\eprf
Finally we prove Lemma~\ref{lem:existenceofh}. In one dimension, one can give an explicit construction as follows. Let $\cos(\eps)<\lam<1$ and set 
\ba
h(t)=\half\sum_{k\not=0}^{\infty}\lam^{|k|}e^{ikt}=\frac{\lam \cos(t)-\lam^{2}}{|1-\lam e^{it}|^{2}}.
\ea
For $t\in \T\setminus[-\eps,\eps]$, clearly this is bounded above by $-\lam(1-\lam)^{-2}(\lam-\cos(\eps))<0$. Thus $\sup_{|t|>\eps}h(t)$ is negative (the precise value is of no importance) while the Fourier coefficients are as desired. We now give a construction valid in any dimension. 
\bprf[Proof of Lemma~\ref{lem:existenceofh}]  Fix $0<\eps<\frac{\pi}{4}$ and let $g,f:\T^{d}\to \R$ be  defined by $g=\one_{(\eps\T)^{d}}$ and $f=(1+\eta)g-\eta$ where $\eta=\frac{\eps^{d}}{1-\eps^{d}}$  is chosen so that $\int_{\T^{d}}fd\lam
=0$. Finally, set $h=f\star f$ be the convolution of $f$ with itself, i.e., $h(t)=\int_{\T^{d}}f(s)f(t-s)d\lam(s)$ (here $\T^{d}$ is treated as a group under addition modulo $2\pi$). 

By choice of $\eta$, we have $\hat{f}(0)=0$. Further, $(g\star g)(t)=\prod_{j=1}^{d}(\eps-\frac{1}{2\pi}|t_{j}|)_{+} $ and $g\star \one=\eps^{d}$. From this it follows that
\ba
(f\star f)(t) &= (1+\eta)^{2}\prod_{j=1}^{d}\l(\eps-\frac{1}{2\pi}|t_{j}|\r)_{+}-2\eta(1+\eta)\eps^{d}+\eta^{2} \\
&=(1+\eta)^{2}\prod_{j=1}^{d}\l(\eps-\frac{1}{2\pi}|t_{j}|\r)_{+}-\eta^{2}
\ea
since $(1+\eta)\eps^{d}=\eta$.

Thus, $h(t)=-\eta^{2}$ for $\|t\|_{\infty}>2\pi\eps$ and $\hat{h}(m)=(\hat{f}(m))^{2}$ which is zero for $m=0$ and non-negative for all $m$. Lastly, $h(0)=(1+\eta)^{2}\eps^{d}$. Dividing $h$ by $h(0)$ gives the desired function.
\eprf


\berk What is it that makes our proof work? Consider two centered Gaussian vectors $X$ and $Y$ in $\R^{n}$ with covariance matrices $\Sigma$ and $\Sigma'$. There are two possible ways to compare $\Sigma$ and $\Sigma'$. Firstly, we may compare them in positive definite order, i.e., $\Sigma\ge \Sigma'$ if $u^{t}\Sigma u\ge u^{t}\Sigma' u$ for all vectors $u$. As the proof of Lemma~\ref{lem:bboundedbelow} shows, in this case, $\Sigma^{-1}\le \Sigma'^{-1}$ and hence,
\ba
\int_{\R_{+}^{n}}\exp\l\{-\half u^{t}\Sigma^{-1}u\r\} du\ge \int_{\R_{+}^{n}}\exp\l\{-\half u^{t}\Sigma'^{-1}u\r\}.
\ea
Although the inequality for the determinant in the denominator of the Gaussian density goes the other way ($\det(\Sigma)^{-\half}\le \det(\Sigma')^{-\half}$), these determinants can be easily bounded (by $\hat{\mu}(0)^{|A|}$ for example) and hence, with a little imprecision, we may  say that if $\Sigma\ge \Sigma'$ in the positive definite order, then   $\P\{X_{i}>0\mb{ for all }i\}$ is smaller than $\P\{Y_{i}>0\mb{ for all }i\}$.   The second comparison  is the one used in Slepian's inequality (entrywise comparison of $\Sigma$ and $\Sigma'$ provided the diagonals are equal). In this case, the much more non-trivial inequality of Slepian gives a comparison of the two probabilities,  $\P\{X_{i}>0\mb{ for all }i\}$ and $\P\{Y_{i}>0\mb{ for all }i\}$. 

These two orderings are rather different from each other, and hence, by mixing them, we are able to compare many more covariance matrices than is possible by either order alone!
\eerk

\section{An example where the gap probability decays faster than exponential}\label{sec:counterexample1}
We give an example to show that the condition of Theorem~\ref{thm:lbdexponential} is necessary. This example is one among a larger class of time series considered by Majumdar and Dhar~\cite{majumdardhar}. Newell and Rosenblatt~\cite{newellrosenblatt} also remark in their paper that if the covariance is not positive, then the gap probability can decay faster than exponential, but they do not give an example. 
\beg \label{eg:derivativeofwhitenoise} Let $d\mu(t)=b(t)dt$ with $b(t)=2-2\cos(t)$. This is the spectral measure of the Gaussian process $X_{n}=\xi_{n}-\xi_{n+1}$ where $\xi_{i}$ are i.i.d. $N(0,1)$. Therefore 
\begin{align*}
H_{X}(n) &= \P\{X_{1}>0,\ldots ,X_{n}>0\}=\P\{\xi_{1}>\xi_{2}>\ldots >\xi_{n+1}\} = \frac{1}{(n+1)!} 
\end{align*}
which decays faster than exponential (to be precise, decays like $e^{-cn\log n}$). Therefore, in general we cannot expect an exponential lower bound.
\eeg
Observe that $b(t)=2-2\cos(t)$ satisfies the hypothesis of the second part of Theorem~\ref{thm:lbdsubexponential} with $p=2$ (and $\eps=\pi$ and $C=10$) and hence, $H_{X}(N)\ge e^{-cN\log N}$. Thus, Theorem~\ref{thm:lbdsubexponential} gives the right lower bound for this example. More generally, consider any finite moving-average process, i.e., a process on $\Z$ of the form $X_{n}=\sum_{k=0}^{r}a_{k}\xi_{n+k}$, where $a_{k}\in \R$ and $\xi_{k}$ is an i.i.d. sequence of standard Gaussians. Then the spectral measure is $b(t)dt$ where $b(t)=\big|\sum_{k=0}^{n}a_{k}e^{ikt}\big|^{2}$ is a trigonometric polynomial. There are two possibilities. 
\benu
\item If $\sum_{k=0}^{r}a_{k}\not=0$, then $b(0)\not=0$ and Theorem~\ref{thm:lbdexponential} applies. We get $H_{X}(N)\ge e^{-cN}$. 
\item If $\sum_{k=0}^{r}a_{k}=0$, then $b(0)=0$, then there is some $p$ such that $b^{(p)}(0)>0$. If $b(t)$ has no zeros in $\T$ other than $0$, then Theorem~\ref{thm:lbdsubexponential} applies and we get $H_{X}(N)\ge e^{-cN\log N}$.\eenu

\section{Proof of Theorem~\ref{thm:lbdsubexponential}}
\label{sec:lbdsubexponential}
To prove Theorem~\ref{thm:lbdsubexponential}, we need a famous result on trigonometric polynomials due to Tur\'{a}n and Remez and Nazarov. For a subset $E\subseteq \T$ and $0< p<\infty$, and $P:\T\to \R$, let $\|P\|_{L^{p}(E)}=(\int_{E}|P(t)|^{p}d\lam(t))^{1/p}$. The following theorem is due to Nazarov~ \cite{nazarovturancomplete}.

\para{Nazarov's complete version of Tur\'{a}n's lemma} There is a number $B$ such that for any trigonometric polynomial in one variable $P(t)=\sum_{k=0}^{n}a_{k}e^{ikt}$, for any measurable subset $E\subseteq \T$ with $\lam(E)\ge \frac{1}{3}$, and for any $0<p\le 2$, we have $\|P\|_{L^{p}(E)}\ge e^{-B(n-1)\lam(\T\setminus E)}\|P\|_{L^{p}(\T)}$. 

\medskip
The original result due to Tur\'{a}n was the inequality $\|P\|_{E}\ge \l(\frac{\lam(E)}{4e}\r)^{n-1}\|P\|_{T}$ for the case when $E$ is an arc in $\T$ (here $\T$ is naturally identified with the circle). The inequality here is for sup-norms while we shall need the comparison of $L^{2}$ norms. Further, Tur\'{a}n's result was valid for arcs only, while Nazarov's inequality is valid for any measurable $E$. Further, in Nazarov's version, as opposed to Tur\'{a}n's original inequality, the bound for ratio between $\|P\|_{E}$ and $\|P\|_{\T}$ goes to $1$ as $\lam(\T\setminus E)\to 0$. All these three features of Nazarov's version of  Tur\'{a}n's inequality are essential to our application below.

\bprf[Proof of Theorem~\ref{thm:lbdsubexponential}] In this proof, let $A=\{1,2,\ldots ,N\}$ and we write $\Sigma_{N}$ for $\Sigma_{A}$ and  $H_{X}(N)$ for $H_{X}(A)$ etc.
\benu
\item Recall from the proof of Lemma~\ref{lem:bboundedbelow} that for $u\in \R^{n}$, with $U(t)=\sum_{k=1}^{N}u_{k}e^{i\<k,t\>}$, we have
\ba
\|u\|^{2} &= \intt_{\T}|U(t)|^{2}d\lam(t) \qquad \mb{ and } \qquad
u^{t}\Sigma_{N}u&\ge \intt_{\T}|U(t)|^{2}b(t)d\lam(t).
\ea
Let $E_{\del}=\{t\in \T\suchthat b(t)\ge \del\}$ so that $\lam(E_{\del})\ge \del$ for some $\del>0$ (since we assume that $b$ is not identically zero). Apply Tur\'{a}n's lemma as stated above to get
\ba
u^{t}\Sigma_{N}u &\ge \del \|U\|_{L^{2}(E_{\del})}^{2}\\
&\ge \del e^{-B(1-\del)(N-1)}\|U\|_{L^{2}(\T)}^{2} \\
&\ge e^{-\gam N}\|u\|^{2}
\ea
where we have absorbed various constants into $\gam$ (hence $\gam$ now depends on $b$). Thus, the smallest eigenvalue of $\Sigma_{N}$ is bounded below by $e^{-\gam N}$. By the third part of Lemma~\ref{lem:bboundedbelow}, we get the lower bound
\ba
H_{X}(N)\ge \l(\frac{1}{2\sqrt{\hat{\mu}(0)}}e^{-\gam N}\r)^{N/2}
\ea
which is at least $e^{-\gam' N^{2}}$ for some $\gam'$. 
\item For the second part, we get a more accurate lower bound for the smallest eigenvalue. For this we again write
\ba
u^{t}\Sigma_{N}u &\ge \del \|U\|_{L^{2}(E_{\del})}^{2}\\
&\ge \del e^{-B\lam(\T\setminus E_{\del})(N-1)}\|U\|_{L^{2}(\T)}^{2} \\
&\ge \del e^{-C\del^{p}N}\|u\|^{2}
\ea
by the assumption that $\lam(\T\setminus E_{\del})\le C\del^{p}$. Choosing $\del=N^{-1/p}$, we get the lower bound $C'N^{-1/p}$ for the smallest eigenvalue of $\Sigma_{A}$. Again invoking the third part of Lemma~\ref{lem:bboundedbelow}, we get the lower bound $H_{X}(N)\ge e^{-\gam N\log N}$ for a constant $\gam$ that depends on $b$ through $C$ and $p$. \qedhere
\eenu
\eprf



\section{An example that (perhaps!) achieves  $e^{-cN^{2}}$}\label{sec:counterexample2}
Let $X$ be the Gaussian process with spectral density $b(t)=1$ for $t\in [\half \pi,\pi]\cup[-\pi,-\half \pi]$ and $b(t)=0$ for $t\in (-\half \pi,\half \pi)$. We have convincing evidence, but not yet a proof, that $H_{X}(N)\le e^{-cN^{2}}$. Note that the covariance kernel in this case is
\ba
K(m)=\begin{cases}
\half & \mb{ if }m=0, \\
0 &\mb{ if }m\mb{ is even and }m\not=0, \\
\frac{1}{\pi m} &\mb{ if }m=3 \; (\mb{mod } 4), \\
-\frac{1}{\pi m} &\mb{ if }m=1 \; (\mb{mod }4).
\end{cases}
\ea
As such we can invert $\Sigma_{N}$ for small $N$ and by numerical experiments on Mathematica for $N\le 24$, we have strong evidence that all entries of $\Sigma_{N}^{-1}$ are positive. Accepting this, it follows that for any $u\in \R_{+}^{N}$, we have $u^{t}\Sigma_{N}^{-1}u\ge \sum_{k=1}^{N}\sig_{N}^{k,k}u_{k}^{2}$, where we use the notation that $\sig_{N}^{i,j}$ is the $(i,j)$ entry of $\Sigma_{N}^{-1}$. Thus,
\ba
H_{X}(N) &\le \frac{1}{\sqrt{\det(\Sigma_{N})}}\prod_{k=1}^{N}\frac{1}{\sqrt{2\pi}}\int_{\R_{+}}e^{-\half \sigma_{N}^{k,k}u^{2}}du \\
&= \frac{1}{2^{N}\sqrt{\det(\Sigma_{N})}\prodd_{k=1}^{N}\sqrt{\sigma_{N}^{k,k}}}.
\ea
Again numerically, we can evaluate the right hand side (call it $\hat{H}(N)$), and it is observed that the points $(k,-\log\hat{H}(k))$ for $1\le k\le 24$, lie very close to the parabola $3.1-0.8x+0.57x^{2}$. This suggests that $H_{X}(N)$ is indeed bounded above by $e^{-cN^{2}}$.

\berk As remarked in the introduction, lower bounds for gap probability in continuous time are stronger than analogous results in discrete time. For upper bounds the reverse is true: ($\P\{X_{t}>0\mb{ for all }t\in [0,T]\}\le \P\{X_{n}>0\mb{ for all }n\le \lfloor T\rfloor\}$). In particular, if it can be proved rigorously for the above process that the gap probability is bounded above by $e^{-cN^{2}}$, then the same holds for the continuous time process $X=(X_{t})_{t\in \R}$ with spectral density (which is now a finite Borel measure on $\R$) $b(t)=\one_{\half \pi\le |t|\le \pi}$. It is worth noting that the process $X$ is not pathological in any sense and in fact it has smooth and even real-analytic sample paths. This is because the covariance function  is real analytic (to see that, either compute the covariance explicitly or use the fact that the Fourier transform of a comapctly supported function is real analytic).
\eerk

\para{Acknowledgement} We thank Satya Majumdar, Deepak Dhar and Naomi Feldheim for  useful discussions on the topics of this paper.

\end{document}